

\documentstyle{amsppt}

\def\a{\kern+.6ex\lower.42ex\hbox{$\scriptstyle \iota$}\kern-1.20ex a}
\define\Cn{\Bbb C^n}
\define\Pn{\Bbb P^n}
\define\Rze{\Bbb R}
\define\Ce{\Bbb C}

\define\dF{\operatorname{d}(F)}
\define\Res{\operatorname{Res}}
\define\mult{\operatorname{mult}}
\NoBlackBoxes

\topmatter
\title
          An inequality for polynomial mappings
\endtitle

\author
            by Arkadiusz P\l{}oski
\endauthor
\abstract
            We give an estimate of the growth
            of a polynonial mapping of $\Cn$.
\endabstract
\address
Department of Mathematics,\newline
Technical University,\newline
Al.~1000\,LPP\,7, 25--314~Kielce, \newline
Poland
\endaddress
\endtopmatter

\document

\head 1.~Main result. \endhead

Let $F=(F_1,\dots,F_n):\Cn \to \Cn$ be a polynomial mapping.
We put \linebreak
$\dF=\text{\#}\,F^{-1}(w)$       
for almost all $w\in \Cn$ and call $\dF$
the geometric degree of ~$F$.

Let $d_i = \deg F_i$\quad for $i=1,\dots,n$. Then
$0 \le \dF \le \prod_{i=1}^nd_i $
if $F_i\not\equiv 0$  for all $i$        (cf.~\cite{6}, p.~434).
Note that $\dF =0$ if and only if the polynomials $F_1,\dots,F_n$ are
algebraically dependent. For any $w\in \Cn$ such that the fiber
$F^{-1}(w)$ is finite we put
$$
\delta_w(F) = \dF - \!\sum_{z\in F^{-1}(w)}\!\!\mult_z(F)
$$
where $\mult_z(F)$ stands for the multiplicity of $F$ at $z$
(cf. \cite{6},~p.~256). We follow the convention that the sum of
an empty family is zero. We have always $\delta_w(F) \ge 0$ for
the finite fibers $F^{-1}(w)$;\quad $\delta_w(F)=0$ if and only if
$F$ is proper at $w$ or $\dF = 0$  (cf. ~\cite{3}).
For any
$z=(z_1,\dots,z_n)\in\Cn$ we put $|z| =\max\limits_{i=1}^n |z_i|$.
We say that an inequality holds for $|z|\gg 1$ if there is a
constant $R>1$ such that it holds for all $z\in\Cn$ such that
$|z| \ge R$.

Let $\widetilde F_i$ be the homogenized polynomial $F_i$. The main
result of this note is

\proclaim{Theorem 1.1}
Suppose that the system of homogeneous equations \newline
$ \widetilde F_1=\dots=\widetilde F_n=0$ has a finite number of
solutions in the projective space $\Pn$.
Then there is a positive constant $C$ such that
$$
|F(z)| \ge C|z|^{-\delta_0(F)} \qquad \text{for } |z|\gg 1.
$$
\endproclaim

The proof of (1.1) will be given at the end of this note. Now, we
will prove the following

\proclaim{Theorem 1.2
{\rm (cf. \cite{1} for the case $n=2$)}}
Assume the assumptions of (1.1) and put
$\mu = \!\sum\limits_{z\in F^{-1}(0)} \!\!\mult_z(F)$.
Then there is a positive constant $C$ such that
$$
|F(z)| \ge C|z|^{\mu -\prod\limits_{i=1}^nd_i +\min\limits_{i=1}^n(d_i)}
\qquad \text{for } |z|\gg 1.
$$
\endproclaim

\demo{Proof of (1.2)}
Let us distinguish two cases.

\item{Case~1.}
$ \dF >\prod\limits_{i=1}^nd_i -\min\limits_{i=1}^n(d_i)$.
Then $F$ is proper by proposition 1.3 of \cite{7}, consequently
$\mu=\dF$ and the inequality follows from theorem 1.10 \cite{7}.

\item{Case~2.}
$\dF \le \prod\limits_{i=1}^nd_i -\min\limits_{i=1}^n(d_i) $.
By definiton of $\delta_0(F)$ we get
$\delta_0(F) \le \prod\limits_{i=1}^nd_i -\min\limits_{i=1}^n(d_i) -\mu $\quad
and the inequality follows from theorem~1.1.
\enddemo

\remark{Remark 1.3}
Let $F=(F_1,\dots,F_n):\Cn\to\Cn$ be a polynomial mapping.
Then the following two conditions are equivalent:

\item{(i)}
the system of equations $F_1=\dots=F_n=0$ has a finite
number of solutions in ~$\Cn$,

\item{(ii)}
there are constants $C>0$ and $q\in\Rze$ such that
$|F(z)| \ge C|z|^q$ \quad for  $|z|\gg 1$
\endremark

\medpagebreak

For $n>2$ the assumption of (1.1) and (1.2) is stronger than (i).
J.~ Koll\`ar  showed in \cite{5} (Prop.~1.10) then we can take
$q=-\prod\limits_{i=1}^nd_i +\min\limits_{i=1}^n(d_i)$\quad in (ii).

\head 2.  Resultant of homogeneous polynomials. \endhead

If $H_1,\dots,H_m$ is a sequence of homogeneous polynomials in
$n+1$ variables, then we denote by $V(H_1,\dots,H_m)$ the set of
all solutions in $\Pn$ of the system $H_1=\dots=H_m=0$.
We will need some properties of the resultant of $n+1$ forms in
$n+1$ variables.

\proclaim{Property 2.1}
If $H_1,\dots,H_{n+1}$ are general forms of degrees
$d_1,\dots,d_{n+1}>0$ in~$n+1$ variables $X=(X_1,\dots,X_{n+1})$,
then their resultant $\Res_X(H_1,\dots,H_{n+1})$ is a polynomial
in coefficients of these forms, homogeneous of degree
$d_1\cdot\dots\cdot d_{n+1}/d_i$ with respect to the
coefficients of $H_i$.
\endproclaim

\proclaim{Property 2.2}
If the coefficients of $H_1,\dots,H_{n+1}$ lie in $\Ce$,\newline
then $\Res_X(H_1,\dots,H_{n+1})=0$\quad if and only if
$V(H_1,\dots,H_{n+1})\neq \emptyset$
\endproclaim

\proclaim{Property~2.3}
Suppose that $H_1,\dots,H_n$ are homogeneous forms of degrees
$d_1,\dots,d_n$                                          
with coefficients in $\Ce$ such that the set
$V=V(H_1,\dots,H_n)\subset\Pn$ is finite.
Let $L$ be a linear form such that $V\cap V(L)=\emptyset$.
For any $p\in V$ we denote by $\mu_p$ the multiplicity of the mapping
$$
\left(\frac{H_1}{L^{d_1}},\dots,\frac{H_n}{L^{d_n}}\right) :
\Pn\setminus V(L) \to \Cn
$$
at $p$.  Then for any homogeneous form $H$ of degree $d>0$ we
have
$$
\Res_X(H_1,\dots,H_n,H) =
   \Res_X(H_1,\dots,H_n,L^d)
   \prod_{p\in V}\left(\frac H{L^d}(p)\right)^{\mu_p}
$$
\endproclaim

The properties (2.1) and (2.2) are well known (cf. ~\cite{4}).
In order to check (2.3) let us assume that the hypersurfaces
$H_i=0$\quad $(1\le i\le n)$ meet transversally i.e. $\mu_p=1$
for all $p\in V$. According to Bezout's theorem $V$ contains
exactly $d_1\cdot\dots\cdot d_n$ points, consequently the product
$\prod\limits_{p\in V}\dfrac H{L^d}(p)$
is a homogeneous form without multiple factors in the
coefficients of $H$. The polynomials
(in the coefficients of the form $H$)\quad
$\Res_X(H_1,\dots,H_n,H)$\quad and
$\prod\limits_{p\in V}\dfrac H{L^d}(p)$
have the same degree equal to $d_1\cdot\dots\cdot d_n$, so by
the~Nullstellensatz there is a constant $R_0\in\Ce$ such that
$\Res_X(H_1,\dots,H_n,H)=R_0\prod\limits_{p\in V}\dfrac H{L^d}(p)$\quad
for every form $H$ of degree $d$.

Putting in the above equality $H=L^d$ we get
$R_0=\Res_X(H_1,\dots,H_n,L^d)$. To check (2.3) in the general case
let us note that the mapping $\Pn\setminus V(L) \to\Cn$ defined in
(2.3) is an analytic branched covering of degree
$d_1\cdot\dots\cdot d_n$. Let $\Omega$ be a Zariski open subset of
$\Cn$ contained in the set of regular values of this mapping. Hence
for any $a=(a_1,\dots,a_n)\in\Omega$ the hypersurfaces
$H_i-a_iL^{d_i}=0$\quad $(1 \le i\le n)$ meet transversally and we
may apply the formula to the homogeneous forms
$H_1-a_1L^{d_1}$,\dots ,$H_n-a_nL^{d_n}$. We obtain the property~2.3
in the general case by passing to the limit when $a\to 0$.

\head 3. Proof of the main result. \endhead

We begin with

\proclaim{Lemma~3.1
{\rm (cf. \cite{2}, lemma~8.2)}}
Let $P(W,T)=P_0(W)T^d+\dots+P_d(W)\in\Ce [W,T]$, %
$d=\deg_TP(W,T)>0$ be a polynomial of $n+1$ variables
$(W,T)=(W_1,\dots,W_n,T)$ such that $P(0,T)\not\equiv 0$.
Let $\delta=d-\deg_TP(0,T)$.
Then there is a positive constant
$C>0$ such that the condition $P(w,t)=0$ implies \quad
$C|t|^{-\delta} \le |w|$ \quad for $|t| \gg1$.
\endproclaim

\demo{Proof of (3.1)}
If $\delta=0$, then the
lemma follows from the theorem on continuity of roots. Let
$\delta>0$. Then $P_0(0)=\dots=P_{\delta-1}(0)=0$ and
$P_{\delta}(0)\neq0$. Suppose $t\neq0$ and put $s=t^{-1}$. The
equation $P(w,t)=0$ can be rewritten in the form
$P_0(w)+\dots+P_{\delta}(w)s^{\delta}+\dots+P_d(w)s^d=0$

By the Weierstrass~Preparation~Theorem we get
$s^{\delta}+Q_1(w)s^{\delta-1}+\dots+Q_{\delta}(w)=0$ near the
origin with holomorphic $Q_j$,\quad $Q_j(0)=0$ for $j=1,\dots,\delta$.
Consequently for small $|s|$, $|w|$ we have:
$|s|^{\delta} \le C_1|w|(|s|^{\delta-1}+\dots+1)\le C_2|w|$
and the lemma follows.
\enddemo

Let $\widetilde Z=(Z_0,Z_1,\dots,Z_n)=(Z_0,Z)$,
$\widetilde F_i=\widetilde F_i(\widetilde Z)$ the homogenized polynomial
$F_i=F_i(Z)$ and let $G=c_1Z_1+\dots+ c_nZ_n$ be a linear form such
that
$$
   V(\widetilde F_1,\dots,\widetilde F_n,Z_0,G)=\emptyset \tag ${*}$
$$
Let $T$ be a variable. We may assume that $d_1,\dots ,d_n>0$. We put
$$
P_G(W,T)=\Res_{\widetilde Z}
   \bigl( \widetilde F_1(\widetilde Z)-W_1Z_0^{d_1},\dots,
          \widetilde F_n(\widetilde Z)-W_nZ_0^{d_n},G(Z)-TZ_0
   \bigr)
$$
{}From the properties of resultant quoted in Section~2 we get immediately

\proclaim{Property 3.2}
$P_G(W,T)\in\Ce [W,T]$,\quad $P_G(w,t)=0$ if and only if the
system of \linebreak
equations $F(Z)=w$, $G(Z)=t$ has a solution in $\Cn$.
In particular $P_G(F(Z),G(Z))=\nomathbreak 0$.
\endproclaim

\proclaim{Property 3.3}
If the fiber $F^{-1}(w)$, $w\in\Cn$ is finite, then           
$\deg_TP_G(w,T)=\!\!\sum\limits_{z\in F^{-1}(w)}\!\!\!\mult_zF$
\endproclaim

Now, we are in a position to prove theorem~1.1. From property~3.3 we
get
$$
  \deg_TP_G(W,T)- \deg_TP_G(0,T)=\delta_0(F) \tag ${*}{*}$
$$
for every form G satisfying (${*}$).

Let $\dF>0$. After a change of coordinates we may assume that the
forms $G=Z_i$\quad $(1\le i\le n)$ satisfy condition (${*}$).
Let $P_i(W,T)=P_G(W,T)$  with $G=Z_i$. 

{}From relations $P_i(F(Z),Z_i)=0$\quad $(1\le i\le n)$, property
(${*}{*}$) and lemma~3.1 we get
$$
  |F(z)| \ge C|z_i|^{-\delta_0(F)}\qquad \text{for } |z_i| \gg 1
$$
for some $C>0$ and (1.1) follows.

Let $\dF=0$. Fix a linear form satisfying (${*}$).
Then $P_G(W,T)=P_0(W)$ and
$F(\Cn)=\{\,w\in\Cn : P_0(w)=0\,\}$ is algebraic. Obviously
$F^{-1}(0)=\emptyset$ (if $\dF=0$ then every fiber of $F$ is empty
or infinite) so there is a $C>0$ such that
$\{\,w\in\Cn : |w|<C\,\} \cap F(\Cn)=\emptyset$, hence
$|F(z)|\ge C$ for all $z\in\Cn$ which proves (1.1) because
$\delta_0(F)=0$ if $\dF=0$.

\enddocument